
\input amstex

\overfullrule=0pt

\documentstyle{amsppt}

 \define\clsp{$\Cal L$ }

\define\emo{$ M_0$}  \define\emosp{$ M_0$ }

\define\bethr{$ B^3$}

  \define\gamsp{$\gamma$ }

\define\ctln{\centerline}

\define\ubr{\underbar}
\define\ssk{\smallskip}
\define\msk{\medskip}

\define\smin{$\setminus$}
 
\define\ap{$\cap$}

\define\up{$\cup$}

\define\sset{$\subseteq$} 

\define\del{$\partial$}
 
\define\delh{$\partial_{h}$}

\topmatter

\title\nofrills Persistent laminations from Seifert surfaces \endtitle
\author  Mark Brittenham$^1$ \endauthor

\leftheadtext\nofrills{Mark Brittenham}
\rightheadtext\nofrills{Persistent laminations from Seifert surfaces}

\affil   University of North Texas\endaffil
\address   Department of Mathematics, Box 305118, University of North
Texas, Denton, TX 76203  \endaddress
\email   britten\@unt.edu \endemail
\thanks \hskip1pc$ ^{1}$ Research supported in part by NSF grants \# 
DMS$-$9400651 and DMS$-$9708411\endthanks

\keywords   essential lamination, Seifert surface, Property 
P\endkeywords

\abstract  We show how an incompressible Seifert surface $F$ 
for a knot $K$ in $S^3$
can be used to create an essential lamination ${\Cal L}_{F}$ in the
complement of each of an infinite class of knots associated to $F$. 
This lamination is persistent for these knots; it remains essential 
under all non-trivial Dehn fillings of the knot complement. This 
implies a very strong form of Property P for each of these 
knots.\endabstract

\endtopmatter

\document

\heading{\S 0 \\ Introduction}\endheading

Essential laminations have found many uses in 3-manifold topology, 
principally in understanding Dehn fillings on 3-manifolds. For 
example, Delman [De] showed that 2-bridge knots have Property P, by 
finding essential laminations in their complements which remain 
essential under all non-trivial Dehn fillings. Such laminations are 
called {\it persistent}, and the knots are
called {\it persistently laminar}. Naimi [Na1],[Na2], independently, proved the same 
result using different methods. In fact, since each of these Dehn filled 
manifolds contains an essential lamination, each has universal cover 
${\Bbb R}^3$ [GO], a property which we have chosen to call 
\ubr{strong} Property P for the knot. Since reducible manifolds do 
not have universal cover ${\Bbb R}^3$, strong Property P implies the 
cabling conjecture for the knot, as well.

More generally, Delman and Roberts [DR], using a mixture of the
above methods, have proved strong 
Property P for non-torus alternating knots. Recently, Wu [Wu] proved 
strong Property P for most algebraic knots, again using essential 
laminations.

In this paper we describe a process for generating knots with strong 
Property P, by building essential laminations in the complement of a 
knot using an incompressible Seifert surface for some other knot. The 
knots which we succeed in doing this for are 
usually not alternating, and are probably not algebraic, since they 
do not seem
to decompose into rational tangles. 

This construction was also discovered, 
around the same time, by Ulrich Oertel (unpublished).

\heading{\S 1 \\ The motivating example}\endheading

The reader is referred to [R] for background in knot theory and Dehn 
surgery. Many of the basic concepts about essential laminations and 
branched surfaces may be found in [GO].

\ssk

The genesis of the constructions we give here is an example due to 
Ulrich Oertel
[Oe1], which was analyzed, from the point of view of tangles, in 
[Br1]. There a branched surface $B$ was 
constructed in the complement of a certain tangle $T_0$, in the 
3-ball \bethr; $B$ carried a lamination \clsp which was essential in 
the complement of any knot obtained by tangle sum with $T_0$ (see 
Figure 1) Further, the lamination remained essential after 
non-trivial Dehn surgery along any of the knots so constructed.

\smallskip

\input epsf.tex

\leavevmode

\epsfxsize=3in
\centerline{{\epsfbox{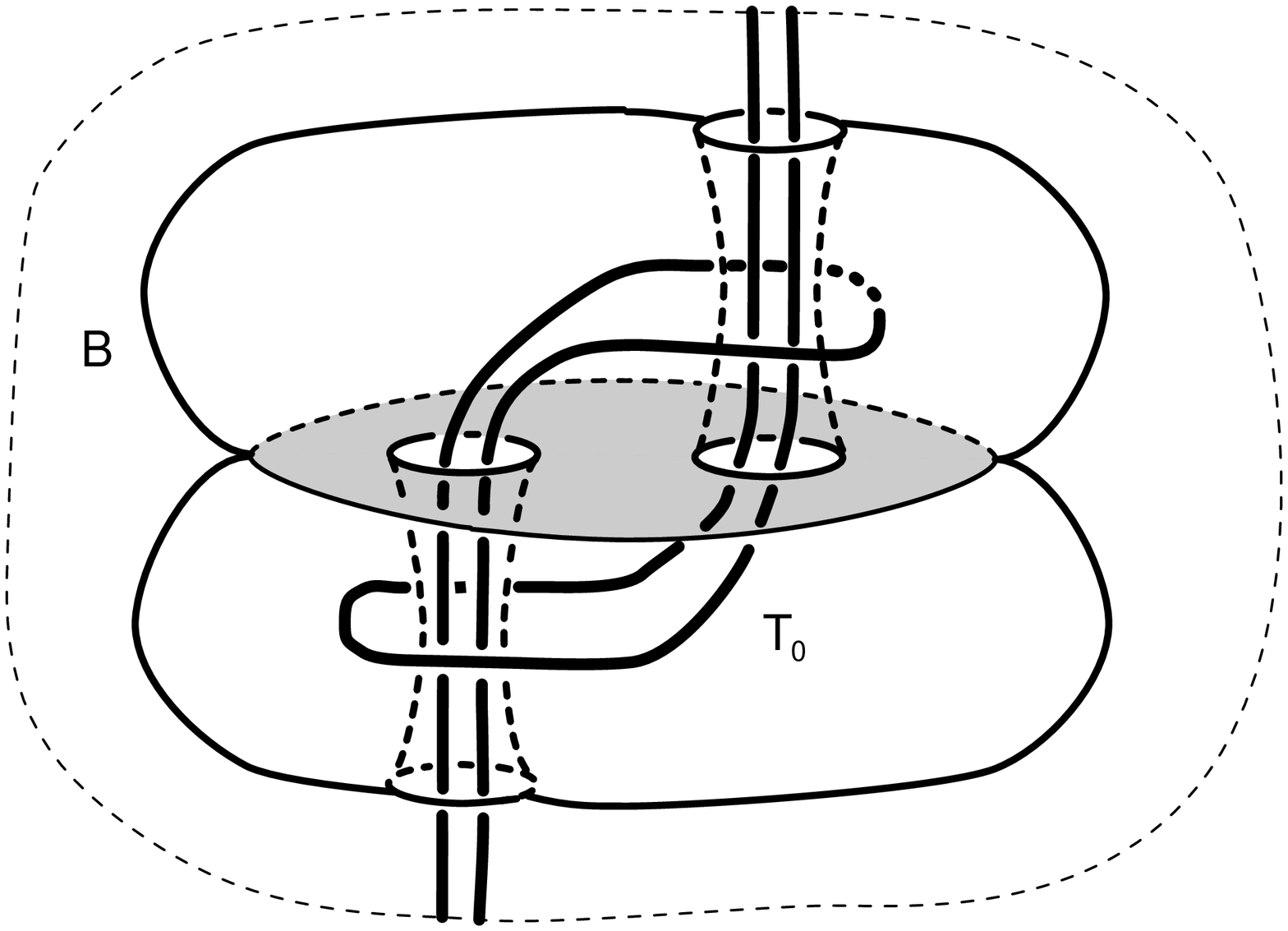}}}

\ctln{Figure 1}

\msk

This lamination and branched surface are, as we shall see, the 
simplest of a long list of laminations and branched surfaces, that 
can be associated to any incompressible Seifert surface for a knot. 
They come, in fact, from applying our construction to a 2-disk 
spanning the unknot. We start with this example, in order to motivate 
the general construction.

\ssk

\leavevmode

\epsfxsize=4.5in
\centerline{{\epsfbox{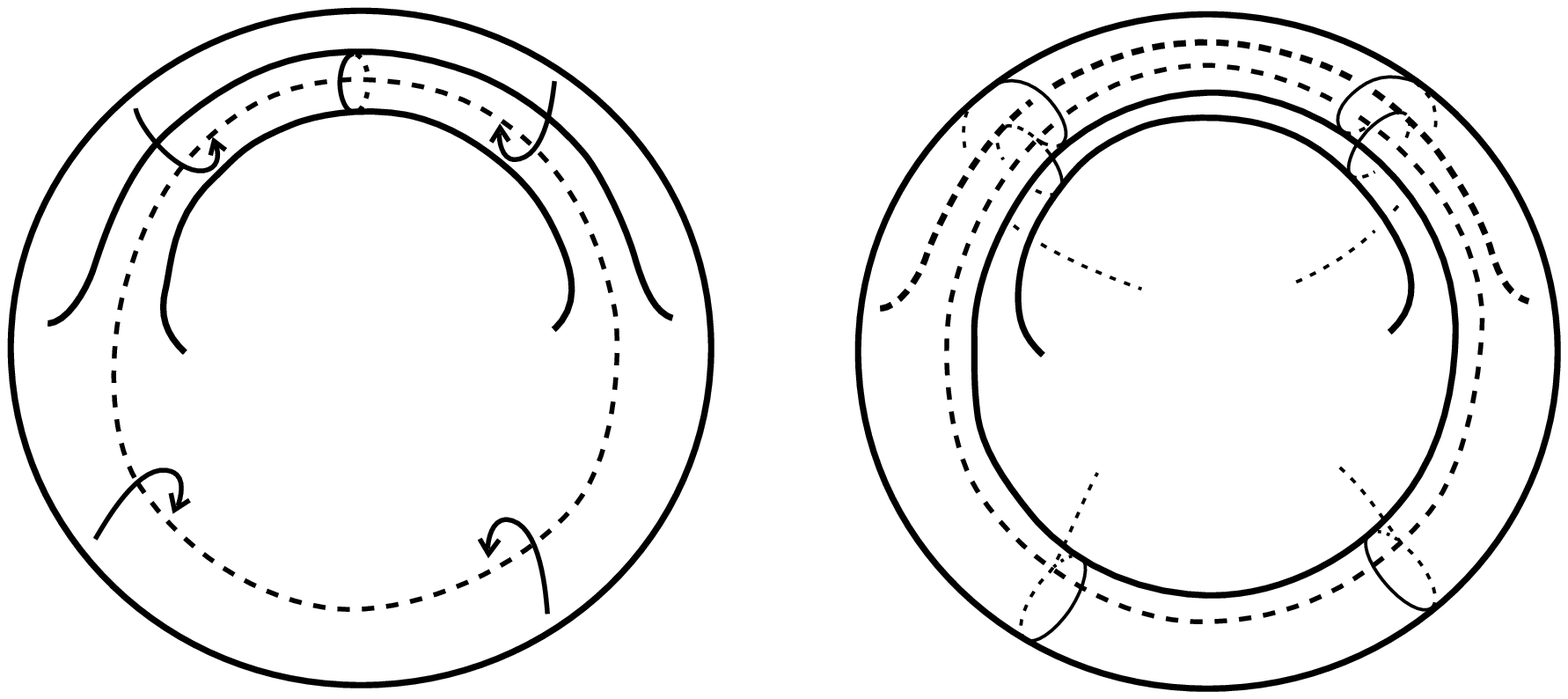}}}

\centerline{Figure 2}

\msk

We can obtain one view of Oertel's branched surface $B$ by 
starting with a 2-disk $D$ spanning the unknot, gluing a tube to it, and 
then gluing the boundary of the resulting punctured torus onto the 
punctured torus, in so doing
creating a branch curve which runs over the tube (see Figure 2). The 
reader is invited to verify that the resulting branched surface is 
isotopic to the one in Figure 1, although this is not really central to what 
follows. If we then string a pair of arcs through both the tube we 
added to $D$ and the tube created when we glued the branched surface 
together (Figure 3), and then complete the pair of arcs to a knot $K$ 
in any way that avoids the branched surface and the two tubes (i.e., 
the two compressing disks for the tubes), then it is not hard to see 
[Br1] that the branched surface $B$ is essential in the complement of 
the knot. Furthermore, since the boundaries of two compressing disks 
of \delh $N(B)$ in $S^3$ are 
not isotopic in \delh $N(B)$ - they live on different components - 
and each compressing disk intersects our knot exactly once, it follows, 
primarily from Menasco's criterion [Me], that our branched surface 
remains essential after any non-trivial Dehn surgery along the knot 
$K$ [Br1]. The knot we have chosen to picture here is in fact the 
twist knot $6_1$, as the reader can verify.

\smallskip

\leavevmode

\epsfxsize=2.5in
\centerline{{\epsfbox{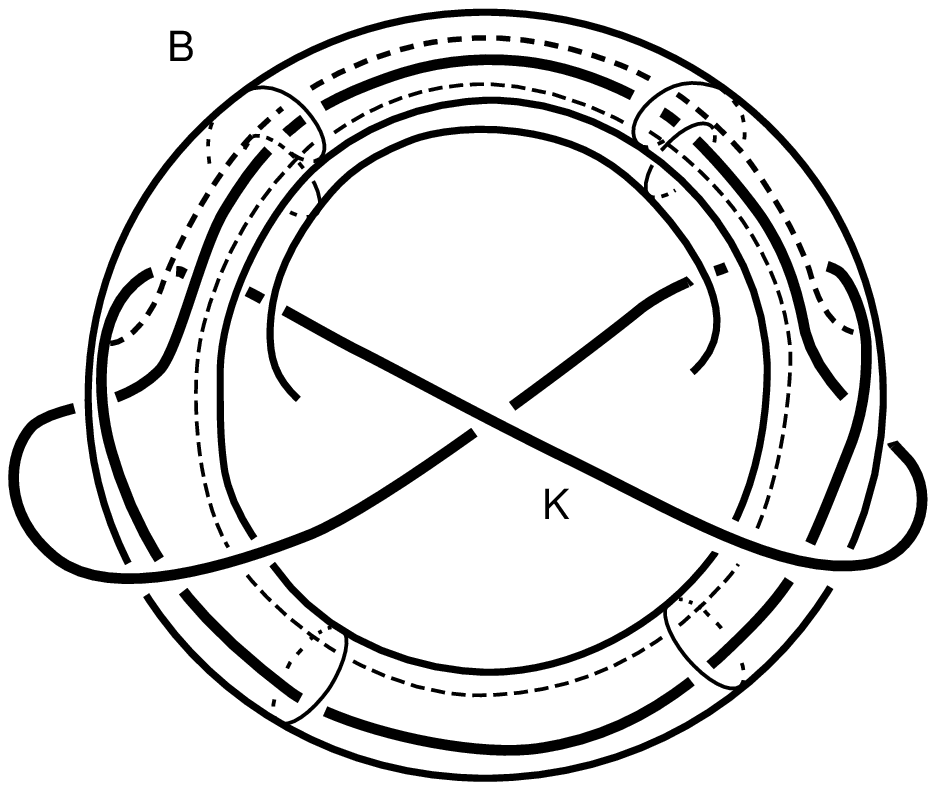}}}

\ctln{Figure 3}

\msk

The fact that we will exploit to find our more general construction 
is that the construction of $B$ did not really use the topology of 
the disk $D$ (or lack of it); it was constructed from $D$ by alterations
taking place
only in a neighborhood of the boundary of $D$. We can therefore apply 
the same
construction to any surface $F$ with boundary in the 3-sphere, i.e., 
to any Seifert surface $F$ for a knot $K$=\del$F$ in the 3-sphere. (The 
construction can be applied to non-orientable surfaces as well, although
conditions which guarantee that the resulting laminations are persistent
are somewhat harder to quantify; see [Br2].)
What we shall see is that if this process is applied to an 
incompressible Seifert surface for a knot, the lamination we create 
is persistent for any knot that we construct in the same manner as above.

\heading{\S 2 \\ The constructions}\endheading

The construction is completely analogous to the one pictured in Figure 
2; we simply forget that we can see the entire disk $D$ and focus on 
a neighborhood of its boundary. Given a Seifert surface $F$ for a knot 
$K$, a neighborhood of its boundary is an annulus; if we attach a 
tube to this annulus, running parallel to half of the knot $K$, and 
then glue the boundary of $F$ (i.e., the knot $K$) to a curve running over the 
tube, and otherwise following the other half of the knot $K$, we get 
a branched surface $B_F$; see Figure 4. There are actually two ways to do this gluing,
to get a branched surface; we choose the one which makes $B_F$ transversely 
orientable, as in the figure. The horizontal boundary of $B_F$ then splits into a positive 
part, which we will call $\partial_+$, and a negative part 
$\partial_-$. This branched surface has a single branch curve $\gamma$, which 
has no self-intersections, and so it is easy to build a `$\beta$-measured' lamination 
\clsp carried by $B_F$ with full support (see [Br1] or [Ro]).

\smallskip

\leavevmode

\epsfxsize=2.5in
\centerline{{\epsfbox{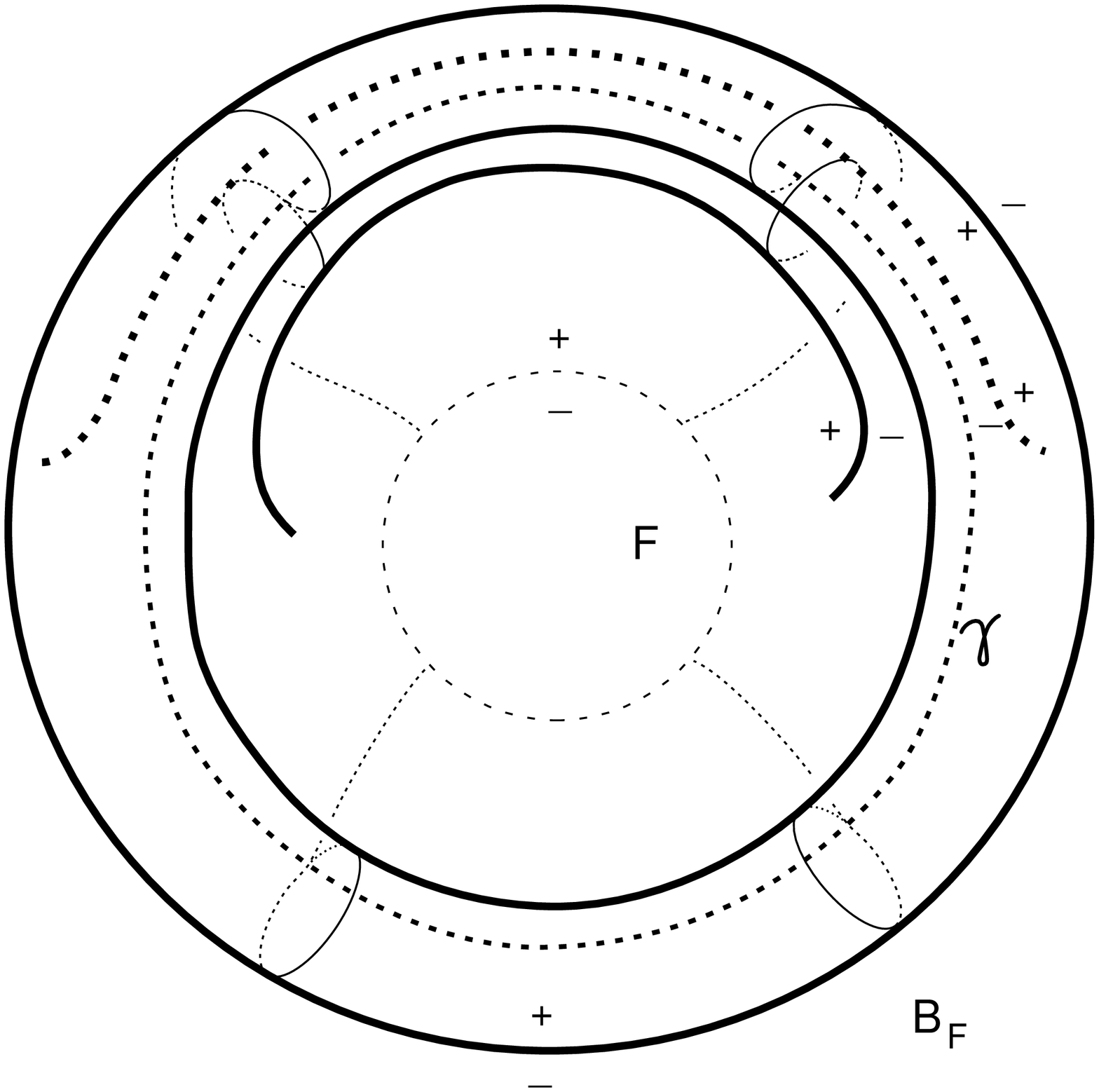}}}

\ctln{Figure 4}

\msk

To understand what will happen next, we need a good picture of what 
the complement of a neighborhood of $B_F$, $M_B$ = 
$S^3\setminus$int$N(B_F)$, looks like. Let $M_K$ = 
$S^3\setminus$int$N(K)$. The idea is to think of the modifications 
described above as taking place in the solid torus neighborhood 
$N(K)$ of $K$ in $S^3$, so that $M_B$\ap$M_K$ = 
$M_K\setminus$int$N(F)$ = $M_F$. $M_F$ can be thought of as a sutured 
manifold, with an annular suture $A$ separating two copies of $F$, 
which we will call $F_+$ and $F_-$. The Seifert surface $F$ is 
incompressible precisely when $F_+$ and $F_-$ are incompressible in 
$M_F$.

\smallskip

\leavevmode

\epsfxsize=4in
\centerline{{\epsfbox{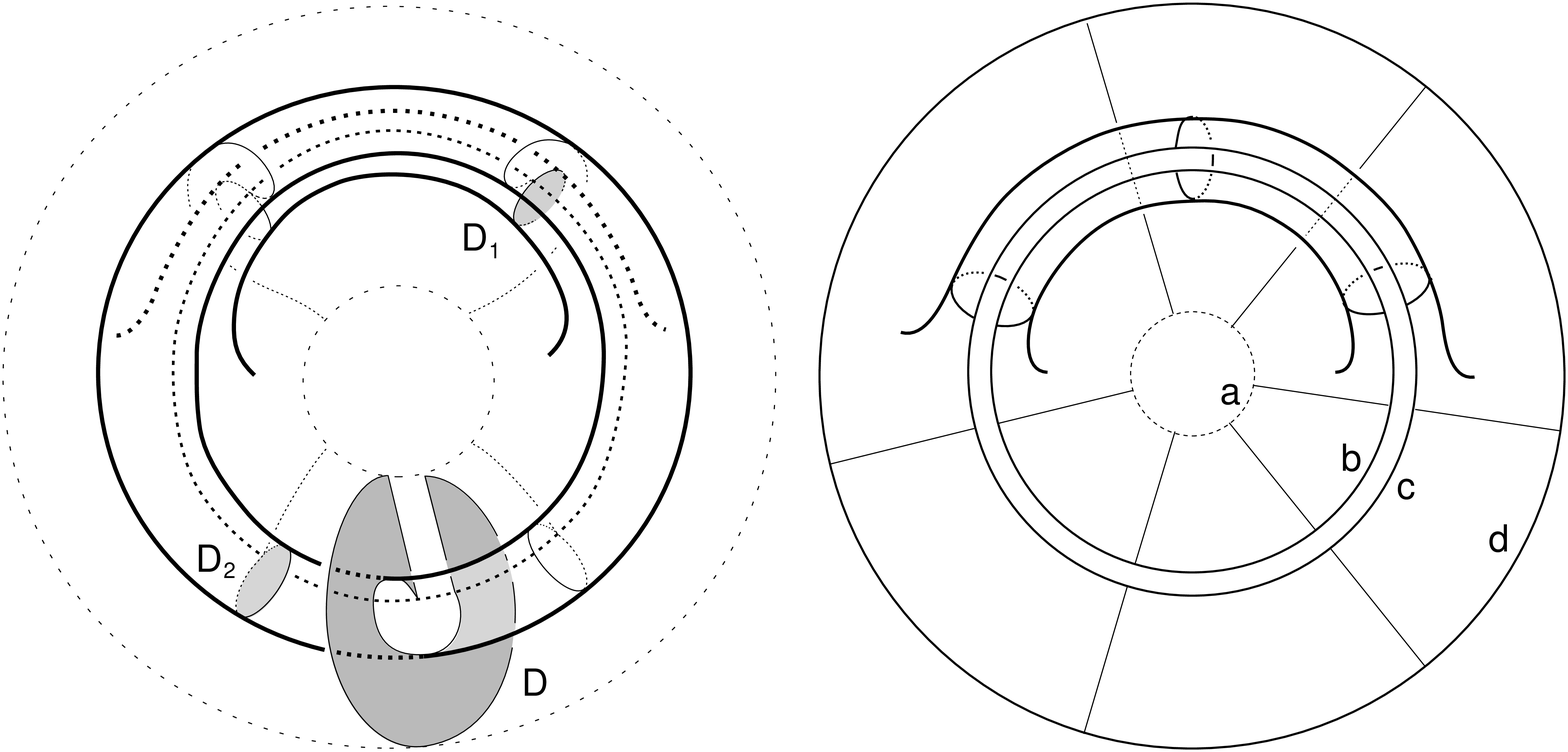}}}

\centerline{Figure 5 \hskip.5in Figure 6}

\msk

To understand $M_B$, it therefore suffices to understand what 
$N(K)$\smin int$N(B)$ = $M_B$\ap$N(K)$ = $N$ looks like. This is a 
sutured manifold, with two sutures, one ($A_1$) being the vertical boundary 
$\partial_v N(B_F)$, and the other ($A_2$) being the annulus 
\del N(K)\smin intN(B) (see Figure 5). By inspection, the complement 
of the two sutures in \del$N$ is a pair of twice-punctured tori; 
they are each built from ($B_F$\ap$N(K)$)\smin$\gamma$, which 
is a 4-punctured sphere (Figure 6), by gluing pairs of boundary circles 
together ($b$ to $c$ for the one from $\partial_+$, $c$ to $d$ for the one from
$\partial_-$). We will call 
these two $\partial$-components $B_+$ and $B_-$. Each of these punctured tori is 
compressible in $N$, via the meridian disks $D_1$,$D_2$ in our two 
tubes (Figure 5). 

Compressing \del$N$ along both of these disks 
yields a new sutured manifold $N_0$, with the same set of sutures, 
whose new positive and negative boundaries are annuli. \del$N_0$ is 
therefore a torus, and so, since $N_0$ is contained in $S^3$ and has 
connected boundary, $N_0$ is irreducible; a reducing
2-sphere must separate boundary components. But \del$N_0$ is still 
compressible; there is a disk $D$ whose boundary meets each suture in 
a single essential arc (Figure 5). Therefore, $N_0$ is a solid torus. 
In particular, because $D$ hits each suture exactly once, $N_0$ is a 
\ubr{product} sutured manifold (annulus)$\times I$. Therefore, $N$ 
is this product sutured manifold with two 1-handles attached, one to the 
positive boundary and one to the negative boundary. Consequently, 

\ssk

$M_B$ $\cong$ $M_F$\up$N$ $\cong$ $M_F$ \up ((annulus)$\times 
I$\up(2 one-handles)) $\cong$ $M_F$\up(2 one-handles)

\msk

\leavevmode

\epsfxsize=2in
\centerline{{\epsfbox{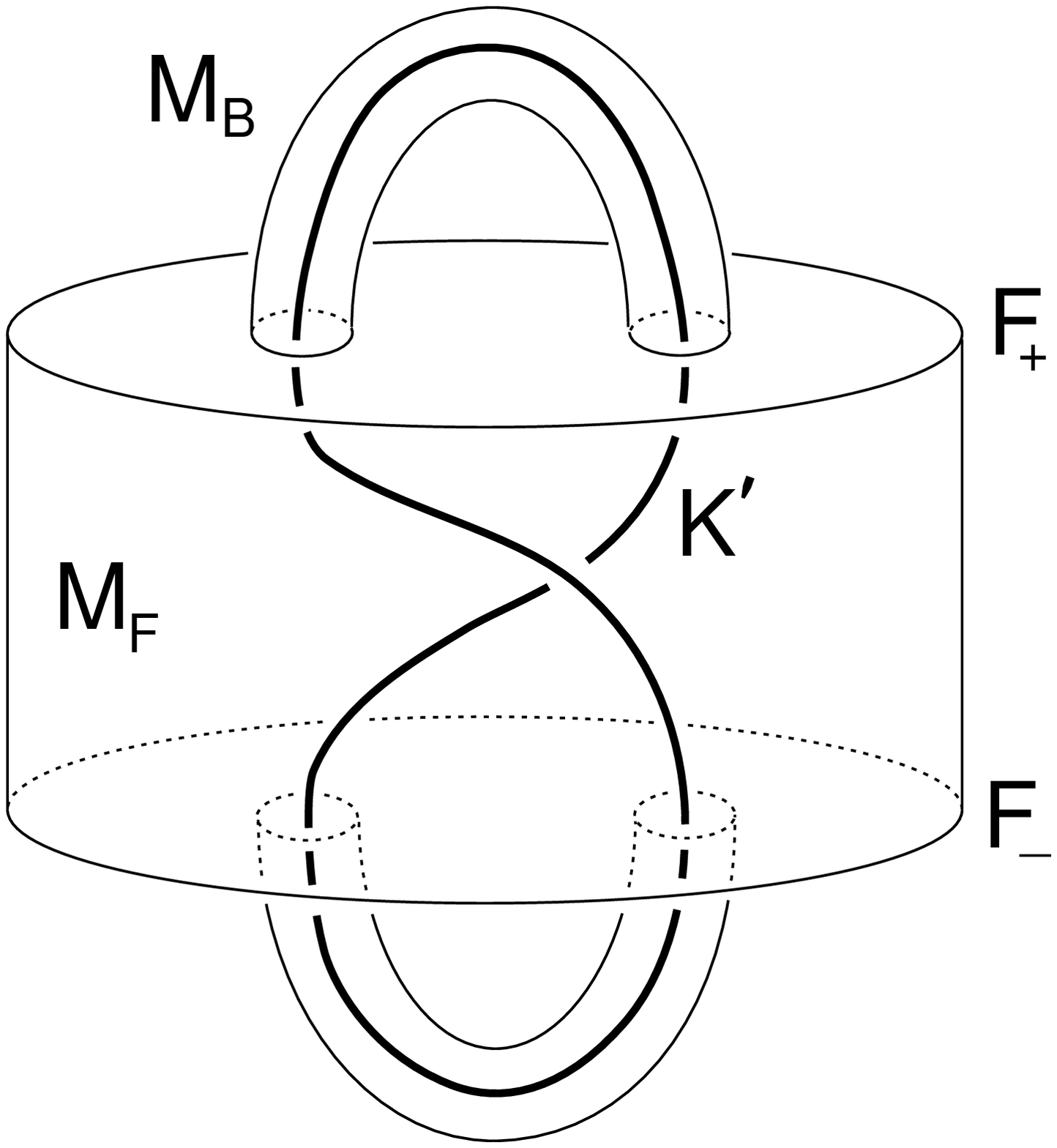}}}

\ctln{Figure 7}

\msk

Therefore, $M_B$ = $M_F$ with two (3-dimensional) one-handles 
attached, one having both ends on $F_+$, the other having both ends 
on $F_-$. \delh N(B) then consists of the surfaces $F_+$ and $F_-$, 
each with a tube attached (increasing their genera by one). It has 
two obvious compressing disks, namely, the meridian disks of the two 
1-handles. So $B_F$ is not essential in $S^3$. To kill these compressing
disks, we will remove a knot $K^\prime$ from 
$M_B$ which intersects each of these disks exactly once (Figure 
7). In other words, we will think of $B_F$ instead as a branched surface in
$S^3$\smin int$N(K^\prime)$=$M_{K^\prime}$. But now we shall easily see 
that if $F$ is an \ubr{incompressible} 
Seifert surface for $K$, then $B_F$ is an essential branched surface 
in $M_{K^\prime}$, and remains essential after every non-trivial Dehn 
filling along $K^\prime$. In other words, $B_F$ is persistent for 
$K^\prime$.

\proclaim{Proposition} 
The branched surface $B_F$ is essential in the complement
of any knot $K^\prime$ build as above.
\endproclaim

Most of the arguments are identical to the proofs from [Br1]. To show 
that $B_F$ is essential in $M_{K^\prime}$, we need to know six things:

\msk

(1)  $B_F$ carries a lamination with full support.

\ssk

This follows, as in [Br1] or [Ro], because the branch curve of $B_F$ does 
not intersect itself.

\msk

(2) $B_F$ does not carry a 2-sphere, and $B$ has no disks of contact.

\ssk

This is because, as in [Br1], $B_F$ has only one sector; $B_F$\smin 
$\gamma$ is connected. Consequently, the branch equations [Oe2] for any
surface carried by $B_F$, or surface of contact for $B_F$, will
be $a=a+a$ or $a=a+a+1$, which have no positive solutions.
In particular, $B_F$ carries no closed surfaces.

\msk

(3) $B_F$ does not carry a compressible torus.

\ssk

This follows for the same reason as (2).

\msk

(4) \emosp = $M_{K^\prime}$\smin int($N(B_F)$) does not have any 
monogons.

\ssk

This follows from the fact that $B_F$ is transversely orientable; the boundary
of a monogon traces out a transverse orientation reversing loop.

\msk

(5) \emosp is irreducible.

\ssk

This is because \del($S^3$\smin int($N(B_F)$)) = \del $N(B_F)$ is 
connected, so $S^3$\smin int($N(B_F)$) is irreducible. Therefore, any 
reducing 2-sphere in \emo\sset $S^3$\smin int($N(B_F)$) bounds a 
3-ball $B^3$ in $S^3$\smin int($N(B_F)$), which therefore contains 
$K^\prime$. But then $K^\prime$ would be null-homotopic in \emo, 
which is impossible since it has non-zero intersection number with 
each of the two compressing disks $D_1$, $D_2$.

\ssk

\leavevmode

\epsfxsize=4.5in
\centerline{{\epsfbox{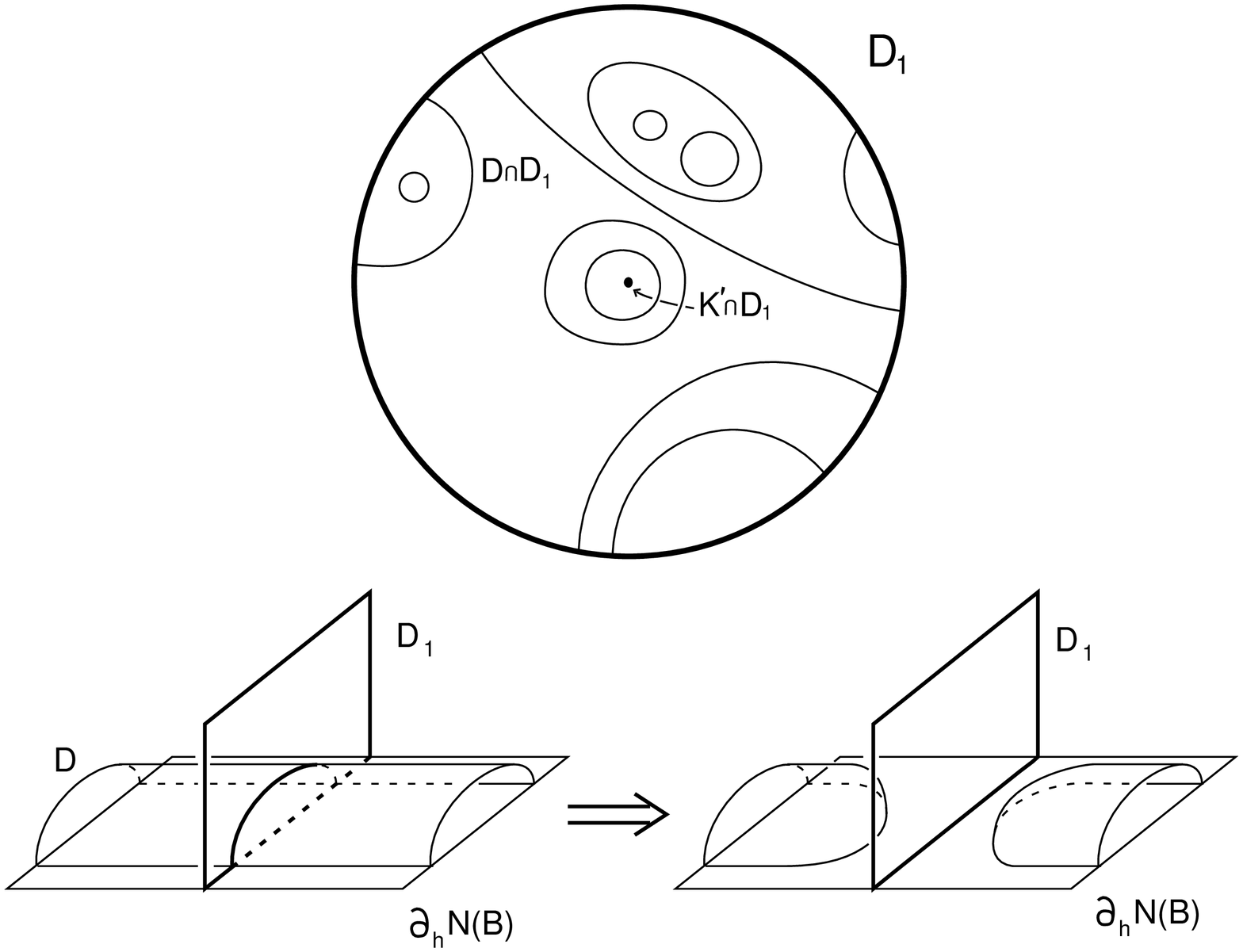}}}

\centerline{Figure 8 \hskip.5in Figure 9}

\msk

(6) The horizontal boundary \delh $N(B_F)$ of $B_F$ is incompressible 
in \emo.

\ssk

Suppose $D$ is a compressing disk for \delh $N(B_F)$ in \emo; in 
particular, $D$\ap$K^\prime$ = $\emptyset$. Consider $D$\ap$D_1$. 
This intersection consists of a finite number of circles and arcs.
The circles come in two types; those which cut off a disk in $D_1$ 
which contains the intersection of $K^\prime$ with $D_1$, and those 
which don't (see Figure 8). We can remove the circles whose disks 
don't meet $K^\prime$ by a standard disk-swapping argument. The arcs 
of intersection can be removed by choosing an outermost arc, and 
using the disk it cuts off from $D_1$ to surger $D$ into two disks 
(see Figure 9). At least one of these is still a compressing disk for 
\delh $N(B_F)$. That disk has fewer arcs of intersection with $D_1$, so 
continuing with that disk will eventually lead to one
with no arcs.

This leaves the circles of intersection which surround the point 
$K^\prime$\ap$D_1$. But the innermost such circle cuts off a disk in 
$D$, which, together with the disk in $D_1$ it cuts off, produces a 
2-sphere which intersects $K^\prime$ exactly once. This is 
impossible, however, in $S^3$.

Therefore, we may assume that $D$ misses $D_1$, and therefore, by 
symmetry, that it misses $D_2$, as well. This means, then that we can 
push $D$ out of the two one-handles that were glued onto $M_F$ to 
create $M_B$. In particular, we may think of \del$D$ as living in 
$F_+$ or $F_-$ (say, $F_+$). But since $F_+$ is incompressible in 
$M_F$, \del$D$ bounds a disk $D^\prime$ in $F_+$. 
In particular, $D\cup D^\prime$ bounds a 3-ball in $M_F$ (since
$M_F$ is irreducible), so $D$ separates $M_F$. If $D^\prime$ 
misses the ends of the one-handle that we glued to $F_+$, then 
$D^\prime$ is contained in \delh $N(B_F)$, a contradiction. One the 
other hand, if $D^\prime$ does contain one or both of these disks, 
then $D$ separates one or both of these disks from both of the disks 
where the other one-handle was glued to $F_-$ (see Figure 10). But 
$K^\prime$ contains a pair of arcs running from the disks on $F_+$ to 
the disks on $F_-$; if they ran from $F_+$ to $F_+$, and from $F_-$ 
to $F_-$, then $K^\prime$ would be a link instead of a knot. So at 
least one of these arcs must intersect $D$, since $D$ is separating. 
But this is also a contradiction.

Therefore, $B_F$ is an essential branched surface in $M_{K^\prime}$. \qed

\smallskip

\leavevmode

\epsfxsize=2.5in
\centerline{{\epsfbox{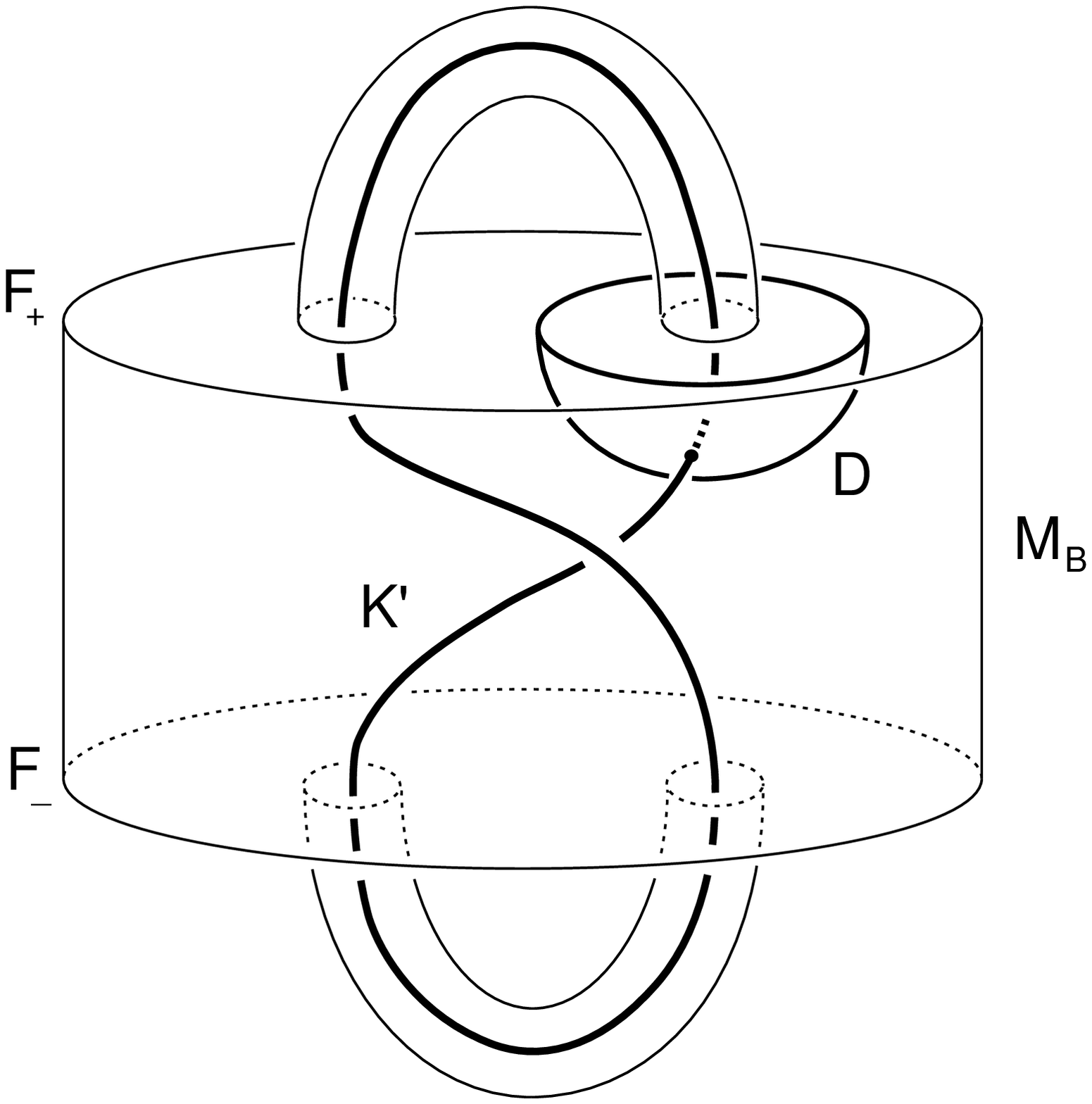}}}

\centerline{Figure 10}

\msk

\heading{\S 3 \\ Persistence}\endheading

The two disks $D_1$ and $D_2$ also allow us to show that $B_F$ 
remains essential after every non-trivial Dehn surgery along 
$K^\prime$. The only properties which do not 
immediately hold after Dehn surgery are (5) and 
(6), because (for (4)) $B_F$ remains transversely orientable, and for all 
of the other properties, $N(B_F)$ has not changed, 
only where it is embedded has, and so our 
previous proofs go through without change. Because \del$D_1$ and 
\del$D_2$ are not isotopic in \delh $N(B_F)$ (they live on distinct 
components), an argument due to Menasco [Me] shows that \delh 
$N(B_F)$ remains incompressible afrer non-trivial Dehn surgery along 
$K^\prime$. It remains to show that $M_{K^\prime}(p/q)$\smin 
int($N(B_F)$) = $M_{(p,q)}$ is irreducible for all $p/q\neq 1/0$. 
Again, this argument is essentially the same as that given in [Br1].

\ssk

Suppose $S$ is a reducing 2-sphere for $M_{(p,q)}$. We may assume 
that $S$ intersects (transversely) the core \gamsp of the solid torus 
glued onto $M_{K^\prime}$, in the fewest number of points. It is then 
standard that $S$\smin int$(N(\gamma))$ = $S^\prime$ is an 
incompressible and \del-incompressible surface in $M_{(p,q)}$\smin 
int$N(\gamma)$ = $M_{K^\prime}$\smin int$N(B_F)$ = $M_2$. The curves 
$S^\prime$\ap\del$N(\gamma)$ are parallel curves of slope $p/q$ on 
\del$N(K^\prime)$\sset\del$M_2$.

The disks $D_1$, $D_2$ meet $M_2$ in annuli $A_1$, $A_2$. Consider 
$S^\prime$\ap$A_i$; it consists of circles and arcs. Trivial circles 
of intersection can be removed by isotopy. The arcs of intersection 
cannot meet the boundary component of $A_i$ coming from \delh 
$N(B_F)$, since $S$ misses \del $M_{(p,q)}$. These arcs of 
intersection are therefore boundary parallel in $A_i$, and so can 
also be removed by isotopy, since $S^\prime$ is \del-incompressible. 
After these isotopies, $S^\prime$ misses the boundary of $A_i$. If 
\del$S^\prime\neq\emptyset$,  then 
\del$S^\prime$\sset\del$N(K^\prime)$ misses a loop in \del$N(\gamma)$ 
which represents a meridian of $K^\prime$, namely $A_i$\ap\del 
$N(\gamma)$, and hence consists of meridional loops. So $p/q$ = 1/0, 
a contradiction.

Therefore \del$S^\prime$=$\emptyset$, in particular, $S^\prime$ = 
$S$\sset$M_2$\sset$M_B$. But since $M_B$ is irreducible, $S$ bounds a 
3-ball in $M_B$. This 3-ball must intersect, hence contain, 
$K^\prime$, since otherwise it is a 3-ball in $M_{(p,q)}$. But this 
implies that $K^\prime$ is null-homotopic in $M_B$, which is 
impossible since it intersects the disk $D_1$ exactly once. So the 
reducing sphere $S$ cannot exist; $M_{(p,q)}$ is irreducible.

\heading{\S 4 \\ Building the knot $K^\prime$, without building
$B_F$}\endheading

\leavevmode

\epsfxsize=4.5in
\centerline{{\epsfbox{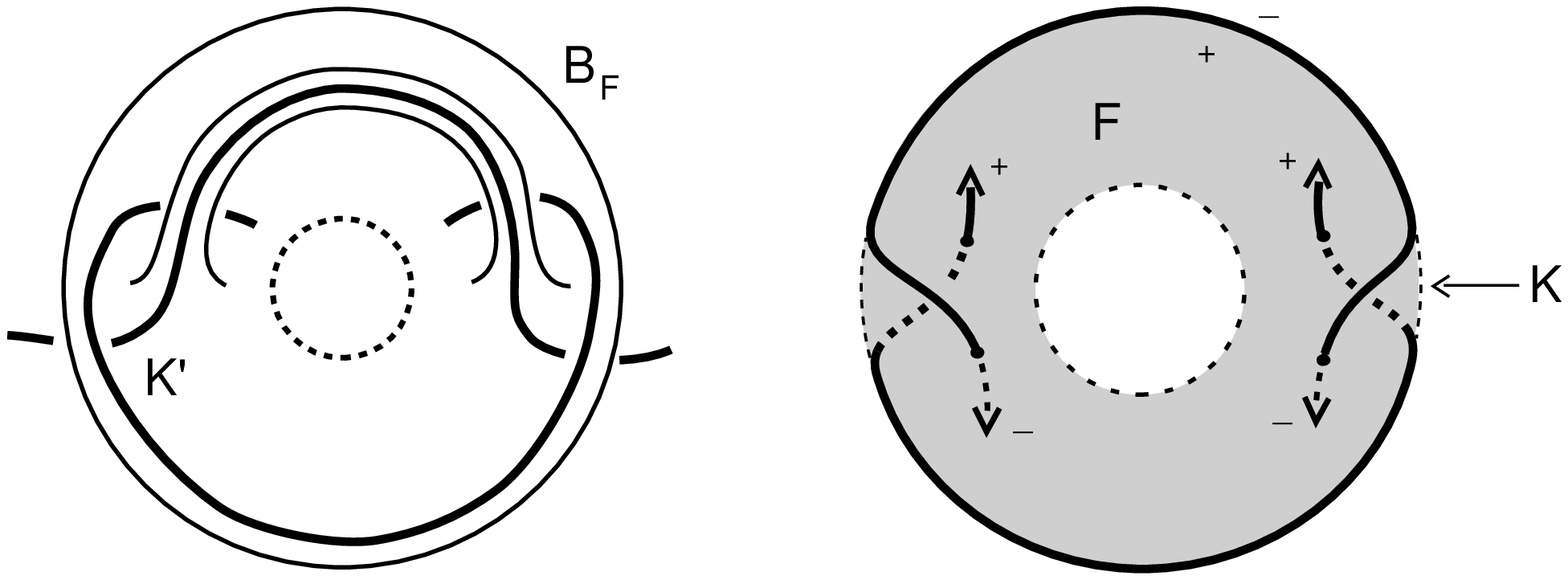}}}

\ctln{Figure 11}

\msk

The knots $K^\prime$ obtained by this construction
can be readily visualized directly from a 
picture of the incompressible Seifert surface $F$ of a knot $K$, 
without constructing $B_F$. The knots $K^\prime$ consist of two arcs 
running through the 1-handles of $S^3$\smin $N(B_F)$, together with 
two arcs in $M_F$ running between these disks on $F_+$ and $F_-$ 
where the two one-handles were attached. Referring to Figures 3 and 11, the 
two arcs in the one-handle can be thought of as a pair of subarcs 
parallel to arcs of $K$, which then pass over and under one another at 
their ends. The ends then reverse direction and pierce the Seifert 
surface $F$ in four points (this reversal is required in order to be 
consistent with the sides of $F$ which the two 1-handles were 
attached to). These two arcs are then completed to a knot in any way 
that misses $F$ (and, technically, neighborhoods of the two added 
crossings). The easiest way to determine which way the two original 
arcs cross each other is to assign a `+' and `$-$' side to $F$, and 
think of one arc remaining slightly on the `$-$' side of $F$, and the 
other remaining slightly on the `+' side of $F$. This gives the 
correct `parity' to the two crossings. Said slightly differently, we 
must make sure that both ends of one arc pierce $F$ from the same side, 
the other ends do so from the opposite side, and arrange the added 
two crossings accordingly. 

The ends of these arcs coming from $K$ 
can then be assigned `+'s and `$-$'s, marking which side of $F$ they 
are emerging from; in order to obtain a knot $K^\prime$, one then 
simply needs to join the four ends by arcs missing $F$, so that the 
+'s are joined to the $-$'s; see Figure 11.

\msk

\heading{\S 5 \\ Some examples}\endheading

The above construction can be carried out for any incompressible Seifert
surface for a knot $K$ . In the case that $K$ is the unknot, it is easy to see 
that all of the knots $K^\prime$ that one can build can be formed by tangle sum 
with a certain algebraic tangle; see [Br1]. For more complicated knots, the
Seifert surfaces will have non-trivial topology, and so their complements will, 
as well. We can take advantage of this extra topology to build a much wider 
variety of persistently laminar knots $K^\prime$.

\smallskip

\leavevmode

\epsfxsize=4.5in
\centerline{{\epsfbox{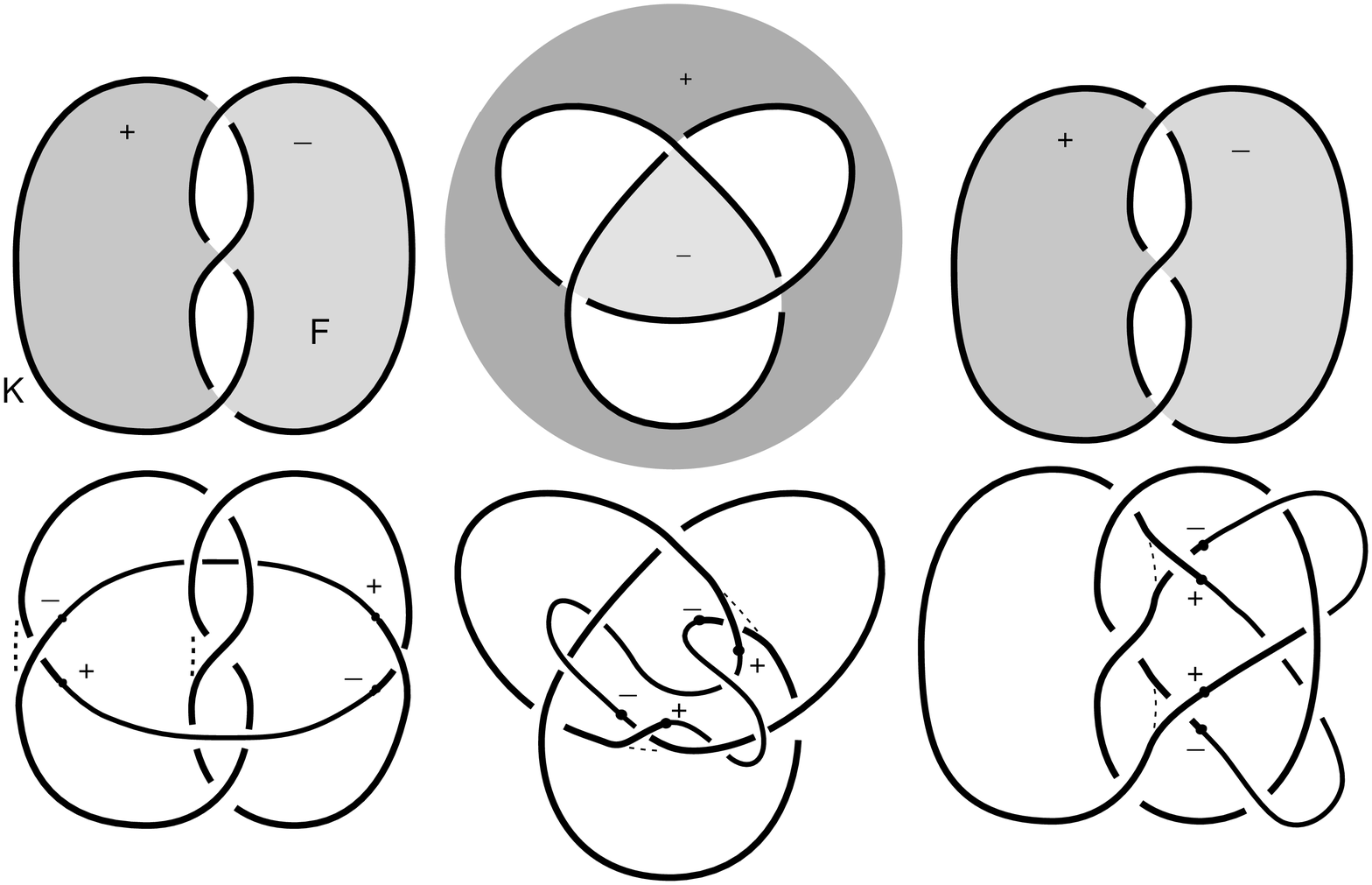}}}

\centerline{Figure 12}

\smallskip

There are (at least) three ways to obtain an incompressible Seifert 
surface for a knot $K$:

\ssk

\noindent (1) run Seifert's algorithm on an alternating projection 
for an alternating knot [Cr];

\ssk

\noindent (2) find a Seifert surface whose genus = 
(span of the Alexander polynomial for the 
knot)/2 [S] - in general, 
genus(K)$\geq$(span)/2;

\ssk

\noindent (3) build a once-punctured torus whose boundary you 
\ubr{know} (for other reasons) is a non-trivial knot.

\msk

Note: all of these in fact give \ubr{least} \ubr{genus} Seifert 
surfaces.

\smallskip

If you have a surface which you suspect is least genus (hence 
incompressible), Gabai's theory of sutured manifold heirarchies [Ga1] 
can, in principle, prove it is. If you have a surface which you 
suspect is incompressible, Haken's normal surface theory (see [JO]) 
can, in principle, prove or disprove it.

\smallskip

We now illustrate this technique for building persistently 
laminar knots with a few examples.

Figure 12 shows several knots that can be built from the Seifert
surface for the trefoil knot. Since the trefoil knot is fibered, 
all of its incompressible Seifert surfaces are isotopic [Th]. The
middle picture shows, however, that different pictures of what 
are really the same surface can be helpful in this construction.
The knots built here turn out to be, from left to right, $9_{46}$, $10_{163}$, 
and $8_9$.

\smallskip

\leavevmode

\epsfxsize=4.9in

\centerline{{\epsfbox{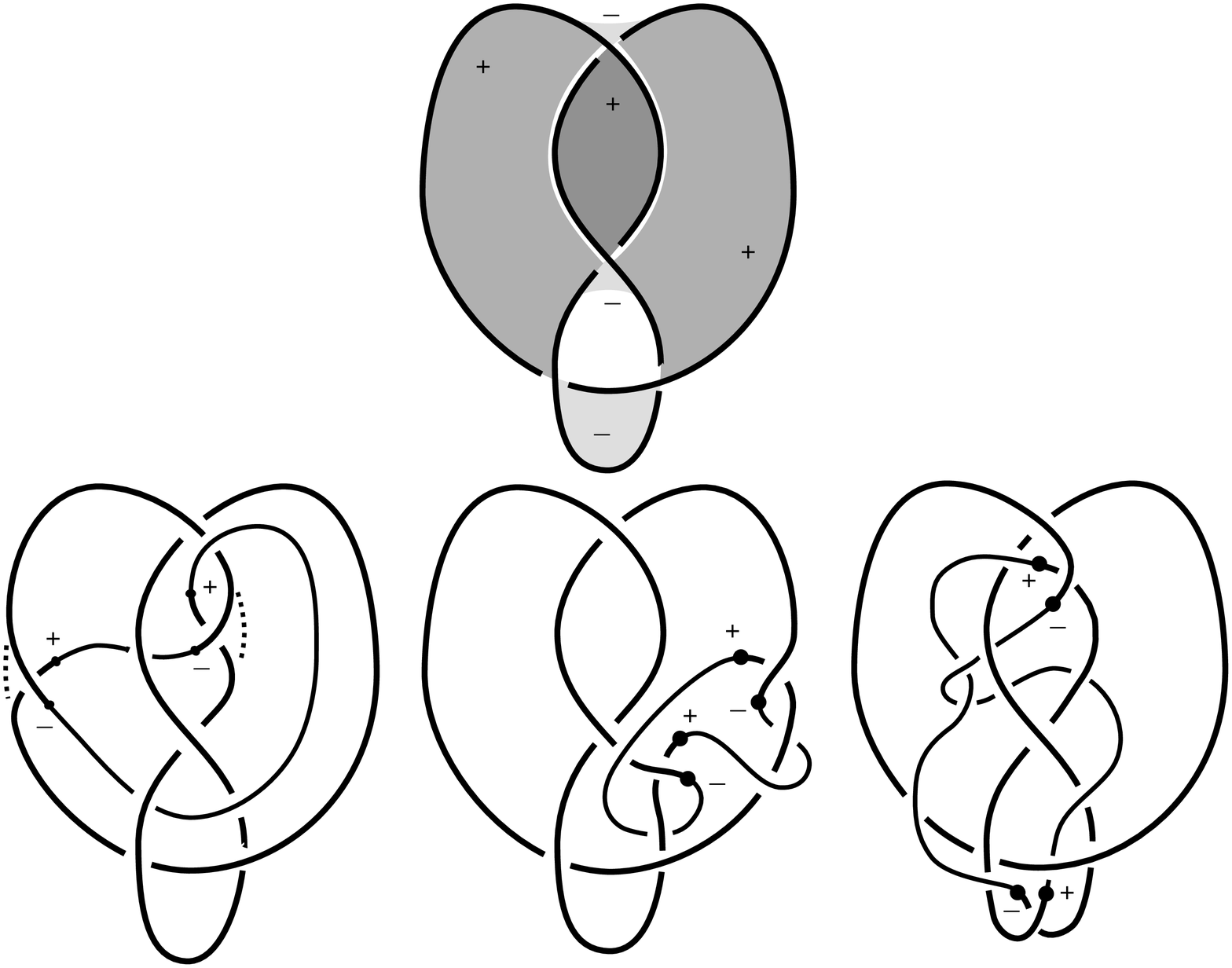}}}

\ctln{Figure 13}

\smallskip

In Figure 13 we build some examples from the figure-8 knot. 
The knots we build are, from left to right, $9_{44}$, $10_{67}$, and $10_{146}$.

Finally, we show how to build an infinite family of knots, by using the 
same, localized, construction. Figure 14 shows how we can alter a knot
in the neighborhood of one of its crossings, if our incompressible Seifert
surface appears as in Figure 14a at the crossing. For example, applying this
construction to the standard Seifert surface for a (2,2$n$+1) torus knot will
build the twist knot with $2n+6$ crossings, i.e., the 2-bridge knots with
continued fraction expansion [2,2$m$] for $m\geq 3$.

\smallskip

\leavevmode

\epsfxsize=3.6in

\centerline{{\epsfbox{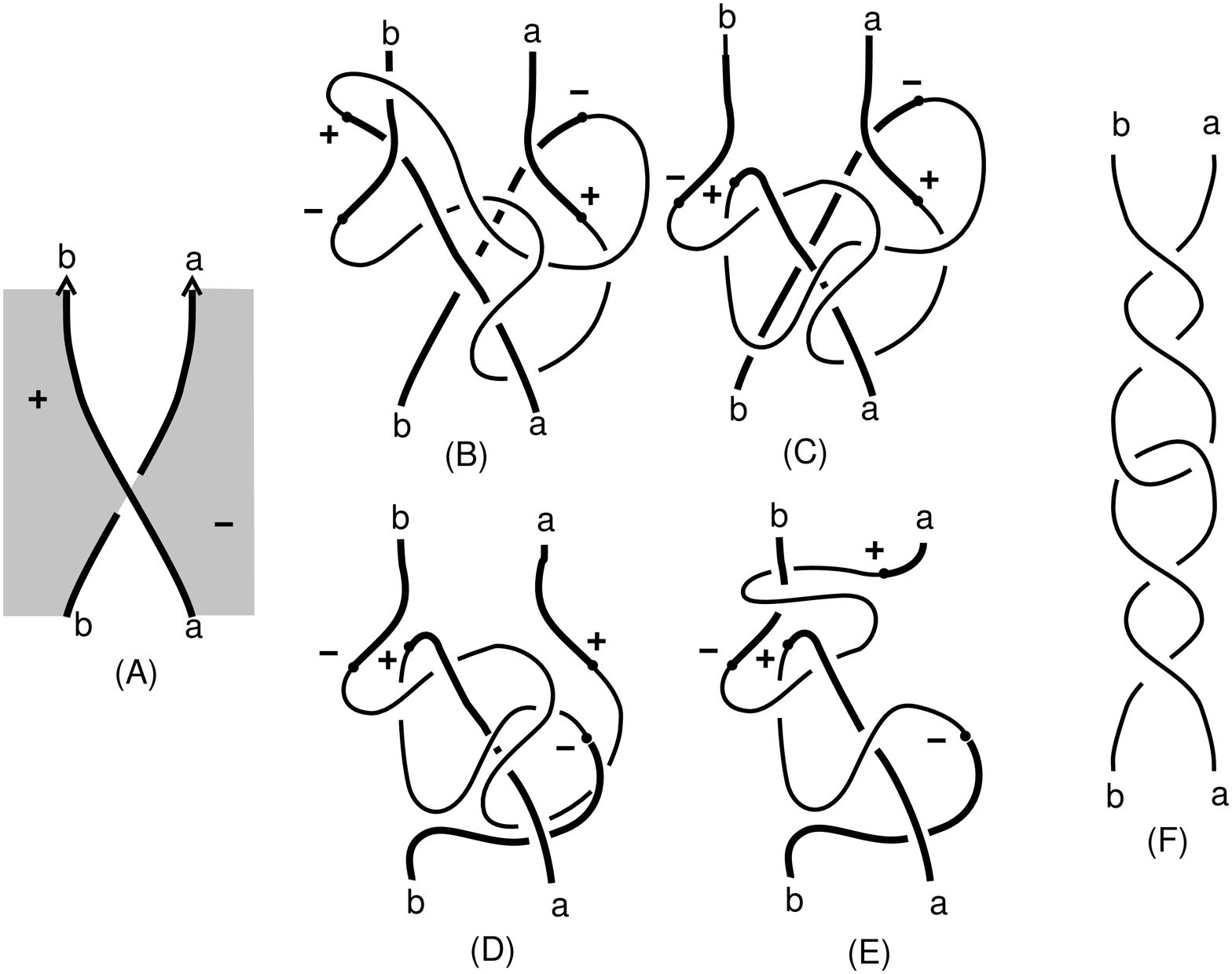}}}

\ctln{Figure 14}

\smallskip

More generally, we can apply this construction to any (2,n) cable of a knot $K$, 
since it is not hard to see that, starting with an incompressible Seifert surface
$F$ for $K$, we can stitch together two parallel copies of $F$, as in Figure 15, to build an 
incompressible Seifert surface $F^\prime$ for the cable. The proof that the resulting 
surface is incompressible comes from the fact that the `holes' in our picture can
be spanned by product disks, in the terminology of sutured manifold theory; any
compressing disk $D$ for $F^\prime$ meets these product disks in circles and arcs. 
The circles can be removed by disk-swapping. (Outermost) arcs, both of whose endpoints 
are on the same side of the product disk, can be removed by surgering $D$ along the 
subdisk each cuts off; one of the two resulting disks must still be a compressing
disk for $F^\prime$. Finally, arcs running from top to bottom cannot exist, since 
the top and bottom 
edges of the product disk are on different sides of $F^\prime$. Once $D$ misses the
product disks, its boundary lives on one of the copies of $F$ that we stitched together;
the incompressibility of $F$ then implies that $\partial D$ bounds a disk in $F$, hence
in $F^\prime$, a contradiction.

\smallskip

\leavevmode

\epsfxsize=4.9in

\centerline{{\epsfbox{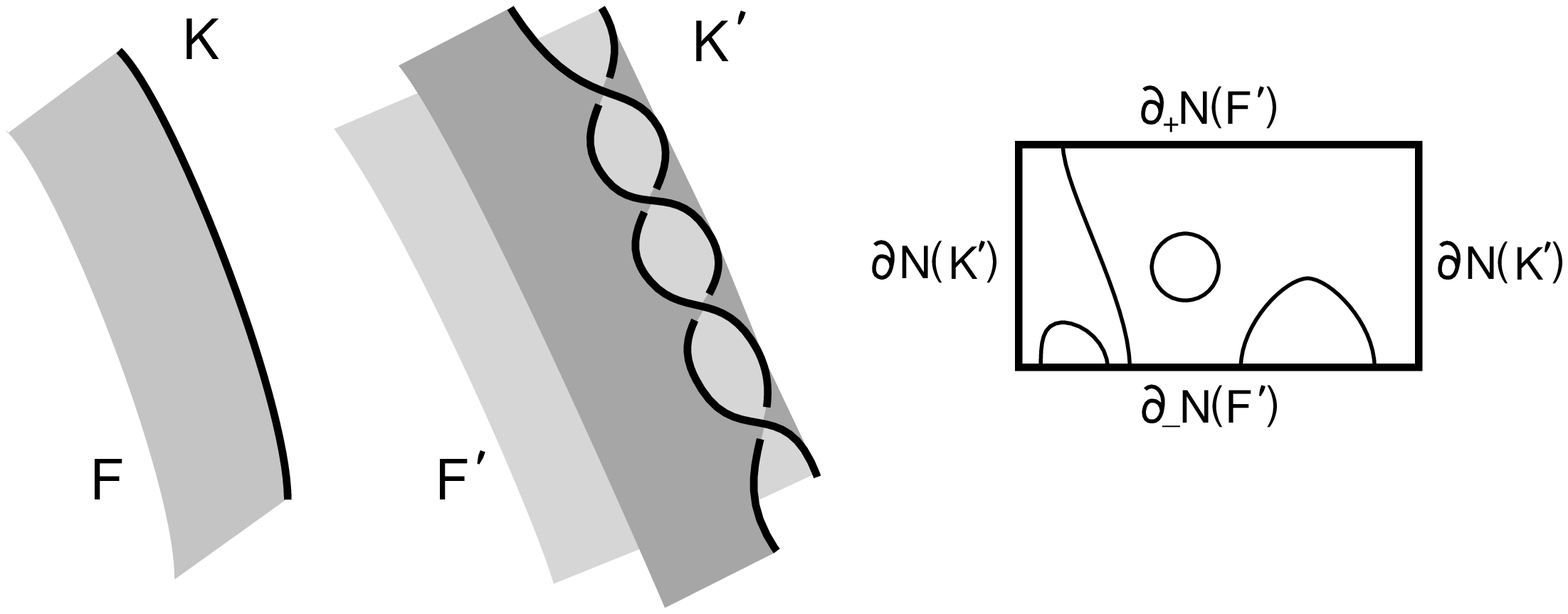}}}

\ctln{Figure 15}

\smallskip

This Seifert surface $F^\prime$ exhibits crossings that we can alter as in Figure 14;
the resulting knot is a twisted Whitehead double of our original knot $K$. Applying this
construction to cables with varying values of $n$, including negative ones (which
yield the mirror image of Figure 14), yields persistent laminations for
slightly less than half of all twisted doubles of $K$. Our twist knots above result 
from applying this construction to a 2-disk spanning the unknot.

\smallskip

Such constructions as above can be carried out in a literally unlimited number of 
ways. From the point of view of [Br1], for example, we can build a wealth of 
persistently laminar tangles by removing any ball in the complement
of $B_F\cup D_1\cup D_2$ meeting our added arcs in a pair of arcs.

In particular, we can continue to add to the list of knots which are known
to have strong Property P.
The table below lists the knots in the standard knot tables which have so far 
been shown to have persistent laminations arising from this construction. A zip file
containing illustrative projections of these knots, in a form readable by
SnapPea [We], can be found at the following web page:

\centerline{http://www.math.unt.edu/$\sim$britten/ldt/knots/knotlst1.html} 

\smallskip

\settabs\+ nnn&nnnnnn&nnnnnn&nnnnnn&nnnnnn&nnnnnn&nnnnnn&nnnnnn&nnnnnn&nnnnnn&nnnnnn&nnnnnn&\cr

\+&
&
&
&
&$ \bold 6_{1}$
&
&
\cr 

\+&
&
&
&
&
&
&
&$ \bold 8_{1}$
&
\cr

\+&
$ \bold 8_{3}$
&
&
&$ \bold 8_{6}$
&$ \bold 8_{7}$
&$ \bold 8_{8}$
&$ \bold 8_{9}$
&$ \bold 8_{10}$ 
&$ \bold 8_{11}$
&
\cr

\+&
$ \bold 8_{13}$
&$ \bold 8_{14}$
&
&$ \bold 8_{16}$
&
&
&
&$ \bold 8_{20}$
&
\cr

\+&
&
&
&
&
&
&$ \bold 9_{7}$
&$ \bold 9_{8}$
&
&
\cr

\+&
&$ \bold 9_{12}$
&
&$ \bold 9_{14}$
&
&
&
&
&$ \bold 9_{19}$
&
\cr

\+&
$ \bold 9_{21}$
&
&
&$ \bold 9_{24}$
&$ \bold 9_{25}$
&$ \bold 9_{26}$
&$ \bold 9_{27}$
&
&
&
\cr

\+&
&$ \bold 9_{32}$
&
&
&
&$ \bold 9_{36}$
&$ \bold 9_{37}$
&
&$ \bold 9_{39}$
&
\cr

\+&
&$ \bold 9_{42}$
&$ \bold 9_{43}$
&$ \bold 9_{44}$
&
&$ \bold 9_{46}$
&$ \bold 9_{47}$
&$ \bold 9_{48}$
&
\cr

\+
&$ \bold 1\bold 0_{1}$
&
&$ \bold 1\bold 0_{3}$
&$ \bold 1\bold 0_{4}$
&$ \bold 1\bold 0_{5}$
&$ \bold 1\bold 0_{6}$
&$ \bold 1\bold 0_{7}$
&
&$ \bold 1\bold 0_{9}$
&$ \bold 1\bold 0_{10}$
\cr

\+
&$ \bold 1\bold 0_{11}$
&$ \bold 1\bold 0_{12}$
&$ \bold 1\bold 0_{13}$
&
&
&$ \bold 1\bold 0_{16}$
&$ \bold 1\bold 0_{17}$
&$ \bold 1\bold 0_{18}$
&
&$ \bold 1\bold 0_{20}$
\cr

\+
&$ \bold 1\bold 0_{21}$
&$ \bold 1\bold 0_{22}$
&$ \bold 1\bold 0_{23}$
&$ \bold 1\bold 0_{24}$
&
&$ \bold 1\bold 0_{26}$
&$ \bold 1\bold 0_{27}$
&
&$ \bold 1\bold 0_{29}$
&$ \bold 1\bold 0_{30}$
\cr

\+
&
&
&
&
&$ \bold 1\bold 0_{35}$
&$ \bold 1\bold 0_{36}$
&
&$ \bold 1\bold 0_{38}$
&
&
\cr

\+
&
&
&$ \bold 1\bold 0_{43}$
&
&
&
&$ \bold 1\bold 0_{47}$
&$ \bold 1\bold 0_{48}$
&
&$ \bold 1\bold 0_{50}$
\cr

\+
&
&
&$ \bold 1\bold 0_{53}$
&
&$ \bold 1\bold 0_{55}$
&
&
&$ \bold 1\bold 0_{58}$
&
&
\cr

\+
&
&$ \bold 1\bold 0_{62}$
&$ \bold 1\bold 0_{63}$
&$ \bold 1\bold 0_{64}$
&$ \bold 1\bold 0_{65}$
&
&$ \bold 1\bold 0_{67}$
&$ \bold 1\bold 0_{68}$
&
&
\cr

\+
&
&
&
&$ \bold 1\bold 0_{74}$
&
&
&$ \bold 1\bold 0_{77}$
&
&$ \bold 1\bold 0_{79}$
&$ \bold 1\bold 0_{80}$
\cr

\+
&
&$ \bold 1\bold 0_{82}$
&$ \bold 1\bold 0_{83}$
&
&
&$ \bold 1\bold 0_{86}$
&$ \bold 1\bold 0_{87}$
&
&
&$ \bold 1\bold 0_{90}$
\cr

\+
&$ \bold 1\bold 0_{91}$
&
&
&$ \bold 1\bold 0_{94}$
&
&
&
&$ \bold 1\bold 0_{98}$
&$ \bold 1\bold 0_{99}$
&
\cr

\+
&
&$ \bold 1\bold 0_{102}$
&
&$ \bold 1\bold 0_{104}$
&
&$ \bold 1\bold 0_{106}$
&
&
&$ \bold 1\bold 0_{109}$
&$ \bold 1\bold 0_{110}$
\cr

\+
&$ \bold 1\bold 0_{111}$
&
&
&
&
&
&
&
&
&
\cr

\+
&
&$ \bold 1\bold 0_{122}$
&
&
&$ \bold 1\bold 0_{125}$
&$ \bold 1\bold 0_{126}$
&$ \bold 1\bold 0_{127}$
&
&
&
\cr

\+
&
&
&$ \bold 1\bold 0_{133}$
&$ \bold 1\bold 0_{134}$
&
&$ \bold 1\bold 0_{136}$
&$ \bold 1\bold 0_{137}$
&
&
&$ \bold 1\bold 0_{140}$
\cr

\+
&$ \bold 1\bold 0_{141}$
&
&$ \bold 1\bold 0_{143}$
&$ \bold 1\bold 0_{144}$
&
&$ \bold 1\bold 0_{146}$
&$ \bold 1\bold 0_{147}$
&$ \bold 1\bold 0_{148}$
&
&$ \bold 1\bold 0_{150}$
\cr

\+
&
&$ \bold 1\bold 0_{152}$
&$ \bold 1\bold 0_{153}$
&$ \bold 1\bold 0_{154}$
&
&
&
&$ \bold 1\bold 0_{158}$
&$ \bold 1\bold 0_{159}$
&
\cr

\+
&
&
&$ \bold 1\bold 0_{163}$
&
&$ \bold 1\bold 0_{165}$
&
\cr

\msk

Combining these results with the results of the constructions listed in the introduction
(see [Ga2]), the list of non-torus knots in the standard knot tables which are not yet
known to have strong Property P becomes remarkably short; 
as of this writing, only the knots $10_{139}$,
and $10_{161}=10_{162}$ (the Perko pair) remain.

\heading{\S 7 \\ Concluding remarks}\endheading

It has long been conjectured that all non-trivial knots in $S^3$ have 
Property P, that is, that non-trivial surgery on the knot can never 
yield a simply-connected manifold. On the other hand, not all knots 
have strong Property P: torus knots, for example, have surgeries with 
finite fundamental group, and cabled knots have reducible surgeries. 
Both of these properties preclude the surgery manifold from having 
universal cover ${\Bbb R}^3$.

The constructions we have described here represent only the simplest
non-trivial sort of branched surface in $S^3$, yet they generate a wealth 
of examples of knots admitting persistent laminations, and hence 
having strong Property P. It is clear that much more can be gained, 
and learned, by applying this construction to other methods of 
building branched surfaces. This has, for example, been carried out 
by Hirasawa and Kobayashi [HK], for some other branched surface constructions.

One question that deserves an answer is whether or not \ubr{every} 
knot with strong Property P admits a persistent lamination. This is 
probably not the case, most likely because 
most exceptional Seifert-fibered spaces do not admit 
essential laminations [Br3],[Cl], although they do have universal cover ${\Bbb 
R}^3$. It is in fact not known that every non-torus alternating knot 
admits a persistent lamination, although as mentioned in the 
introduction, they do have strong Property P. There is at least one 
way to show that a knot does not admit a persistent lamination [BNR], 
but it requires knowing that one of the surgery manifolds is the 
`right' kind of graph manifold. Whether or not this is the case for 
any alternating knot might be an interesting topic of study.

Another question which we \ubr{can} answer is whether or not, for all 
of the knots we can construct by this method, the surgery manifolds 
are in fact all Haken, i.e., contain an incompressible closed 
surface. This construction can easily generate a 2-bridge 
knot from certain other 2-bridge knots; for example, our last construction 
generated twist knots from 2-bridge torus knots. Since a 
2-bridge knot exterior contains no closed, non-\del-parallel 
incompressible surfaces [HT], it follows [Ha] that all but finitely 
many Dehn surgeries on these knots are non-Haken. 

\smallskip

\centerline{\bf Acknowledgements}

\smallskip

The author wishes to thank Rachel Roberts and Ramin Naimi for several helpful 
conversations during the preparation of this work.

\enddocument